\documentclass[12pt]{article}
\usepackage{amsmath,amssymb}
\def\R{\mathbb{R}}
\def\Z{\mathbb{Z}}
\def\arcsinh{\mathrm{arcsinh}}
\def\arctanh{\mathrm{arctanh}}
\newtheorem{theorem}{\hspace*{\parindent}Theorem}
\newtheorem{lemma}{\hspace*{\parindent}Lemma}
\newtheorem{corollary}{\hspace*{\parindent}Corollary}

\title{Parameter convexity and concavity of generalized trigonometric functions}
\author{D.B.\:Karp$^{\rm a,b}$\footnote{Corresponding author. E-mail: D. Karp -- \emph{dimkrp@gmail.com}, E.\:Prilepkina --  \emph{pril-elena@yandex.ru}}~~and E.G.\:Prilepkina$^{\rm a,b}$
\\[10pt]\small{\textit{$\phantom{1}^a$Far Eastern Federal University, 8 Sukhanova street, Vladivostok, 690950, Russia}}\\\small{\textit{$\phantom{1}^b$Institute of Applied Mathematics, FEBRAS, 7 Radio Street, Vladivostok,  690041, Russia}}}
\date{}
\def\R{\mathbb{R}}
\def\Z{\mathbb{Z}}
\def\arcsinh{\mathrm{arcsinh}}
\def\arctanh{\mathrm{arctanh}}
\begin{document}
\maketitle

\begin{center}
\parbox{12cm}{
\small\textbf{Abstract.} We study the convexity properties of the generalized trigonometric functions considered as functions of parameter.  We show that $p\to\sin_p(y)$ and $p\to\cos_p(y)$ are log-concave on the appropriate intervals while $p\to\tan_p(y)$ is log-convex.  We also prove similar facts about  the generalized hyperbolic functions.  In particular, our results settle the major part of a conjecture put forward in \cite{BBV}.}
\end{center}

\bigskip


\bigskip


\bigskip

\paragraph{1. Introduction and preliminaries.}
 The symbol $\R_+$ will mean $[0,\infty)$.  There are several ways in the literature to define generalized
trigonometric functions (see, for instance,
\cite{EdmGurkaLang-JAT,EdmLang,Lindqvist95,Lind-Peetre2003,Shelupsky}).
We will stick with the definition adopted in the book
\cite{LangEdm-book}.  For $p>0$ define a differentiable function
$F_p:[0,1)\to\R_+$ by
\begin{equation}\label{eq:arcsin}
F_p(x)=\int_0^x(1-t^p)^{-1/p}dt.
\end{equation}
Clearly, $F_2=\arcsin$ so that $F_p$ can be viewed as generalized arcsine $F_p(x)=\arcsin_p(x)$. Since $F_p$ is strictly increasing it has an inverse denoted by $\sin_p$.  In all the references we could find the range of $p$ is restricted to $(1,\infty)$ because only in this case $\sin_p(x)$ can be made periodic like usual sine. Nothing prohibits, however,  defining $\sin_p(x)$ for all $p>0$, so we will be dealing with such generalized case here. If $p>1$ the function $\sin_p(x)$  is defined on the interval $[0,\pi_p/2]$, where
$$
\pi_p=2\int_0^1(1-t^p)^{-1/p}dt=\frac{2\pi}{p\sin(\pi/p)}.
$$
It is convenient to extend the above definition by setting $\pi_p=+\infty$ for $0<p\leq{1}$. We will adopt this convention throughout the paper.  In this way the function $y\to\sin_p(y)$ is strictly increasing on $[0,\pi_p/2]$ with $\sin_p(0)=0$ and $\sin_p(\pi_p/2)=1$ in analogy with the usual sine.  It is easily seen that $p\to \pi_p$ is strictly decreasing on $(1,\infty)$ and maps this interval onto itself.  For $p>1$ the definition is extended to $[0,\pi_p]$ by
$$
\sin_p(y)=\sin(\pi_p-y)~\text{for}~y\in[\pi_p/2,\pi_p];
$$
further extension to $[-\pi_p,\pi_p]$ is made by oddness; finally $\sin_p$ is extended to the whole $\R$ by $2\pi_p$ periodicity. If $p\in(0,1]$ the inverse to $F_p(x)$ from (\ref{eq:arcsin}) is defined on $\R_+$ and we just need oddness to extend the definition to the whole real line.  The limit cases are (see also \cite{BushEdm-Rocky}):
\begin{equation}\label{eq:sin-limits}
\sin_0(y)=0~\text{on}~\R,~~~\sin_{1}(y)=1-e^{-y}~\text{on}~\R_+,~~~\sin_{\infty}(y)=y~\text{on}~[0,1].
\end{equation}
Since
\begin{equation}\label{eq:sineprime}
\frac{d}{dy}\sin_p(y)=\left(\frac{dF_p(x)}{dx}\right)^{-1}_{|x=\sin_p(y)}=\left(1-[\sin_p(y)]^p\right)^{1/p}
\end{equation}
we get $\sin_p'(0)=1$ and $\sin_p'(\pi_p/2)=0$ which shows that $\sin_p(y)$ is continuously differentiable on $\R$ for all $p>0$.   The continuous derivative above is naturally called the generalized cosine:
\begin{equation}\label{eq:cosp}
\cos_p(y):=\frac{d}{dy}\sin_p(y),~~y\in\R.
\end{equation}
When $y\in[0,\pi_p/2]$ (for $p>1$) and $y\in\R_+$ (for $0<p\leq1$) we can also define $\cos_p(y)$ by the right hand side of (\ref{eq:sineprime}) which leads to an integral representation for $\arccos_p$:
$$
\cos_p(y)=x=\left(1-[\sin_p(y)]^p\right)^{1/p}~\Rightarrow~y=\arcsin_p\left((1-x^p)^{1/p}\right)=\!\!\!\!\!\!\int\limits_0^{(1-x^p)^{\frac{1}{p}}}\!\!\!\!\!\frac{dt}{(1-t^p)^{1/p}},
$$
or, by substitution $s=(1-t^p)^{1/p}$,
\begin{equation}\label{eq:arccosp}
y=\arccos_p(x)=\!\!\!\!\!\!\int\limits_0^{(1-x^p)^{\frac{1}{p}}}\!\!\!\!\!\frac{dt}{(1-t^p)^{1/p}}=\int\limits_{x}^{1}\frac{s^{p-2}ds}{(1-s^p)^{1-\frac{1}{p}}},~~0\le{x}\le1.
\end{equation}
The function $\cos_p$ can now be defined on $[0,\pi_p/2]$ as the inverse function to $\arccos_p$ and extended to $\R$ by evenness and $2\pi_p$ periodicity.  Limiting values for $p=0,1,\infty$ can be obtained by differentiating (\ref{eq:sin-limits}). Pursuing an  analogy with trigonometric functions further, the generalized tangent function is defined by
\begin{equation}\label{eq:tanp}
\tan_p(y)=\frac{\sin_p(y)}{\cos_p(y)},
\end{equation}
where $y\in\R\!\setminus\!\!\{(\Z+1/2)\pi_p\}$ if $p>1$.  If $0<p\leq1$  the function $\tan_p(y)$ is continuous on $\R$.
It is easy to show by differentiation that $\tan_p(y)$ is  the inverse function to
\begin{equation}\label{eq:arctanp}
\arctan_p(x)=\int_{0}^{x}\frac{dt}{1+t^p},~~0\le{x}<\infty,
\end{equation}
extended from $[0,\pi_p/2]$  to $[-\pi_p/2,\pi_p/2]$ by oddity and further by periodicity \cite{EdmLang,LangEdm-book}.
The limit cases are computed to be:
$$
\tan_0(y)=2y~\text{on}~\R,~~~\tan_{1}(y)=e^{y}-1~\text{on}~\R_+,~~~\tan_{\infty}(y)=y~\text{on}~[0,1].
$$

In a similar fashion one can define the hyperbolic sine $\sinh_p(y)$ on $\R_+$ as the inverse function to the integral ($p>0$)
$$
y=\arcsinh_p(x)=\int\limits_0^x\frac{dt}{(1+t^p)^{1/p}}, ~~~0\le{x}<\infty.
$$
This definition is extended to negative values of $y$ by $\sinh_p(y)=-\sinh_p(-y)$. Here the value $p=1$ does not represent any additional difficulties.   Further, we can define the hyperbolic cosine by
$$
\cosh_p(y)=\frac{d}{dy}\sinh_p(y)=(1+[\sinh_p(y)]^p)^{1/p},~~~~0\le{y}<\infty,
$$
which is extended to negative values of $y$ by $\cosh_p(y)=\cosh_p(-y)$.  This leads to the following identity
\begin{equation}\label{eq:coshsinh}
[\cosh_p(y)]^p-|\sinh_p(y)|^p=1,~~~~~y\in\R.
\end{equation}
Finally, the hyperbolic tangent is naturally defined by
$$
\tanh_p(y)=\frac{\sinh_p(y)}{\cosh_p(y)},~~~~0\le{y}<\infty.
$$
Differentiating (\ref{eq:coshsinh}) with respect to $y$ we derive
$$
\frac{d}{dy}\cosh_p(y)=\frac{[\sinh_p(y)]^{p-1}}{[\cosh_p(y)]^{p-2}}=\frac{[\sinh_p(y)]^{p-1}}{(1+[\sinh_p(y)]^p)^{(p-2)/p}},
~~~y\in\R_+,
$$
and differentiating the definition of $\tanh_p(y)$ with respect to $y$ we get after simplifications
$$
\frac{d}{dy}\tanh_p(y)=\frac{1}{[\cosh_p(y)]^p}
$$
implying
$$
\frac{d}{dx}\arctanh_p(x)=\frac{1}{1-x^p}.
$$
Hence, $\tanh_p(y)$ can be alternatively defined on $\R_+$ as the inverse function to the integral
\begin{equation}\label{eq:arctanh}
y=\arctanh_p(x)=\int_{0}^{x}\frac{dt}{1-t^p},~~~0\le{x}<1
\end{equation}
 - the  definition we will use below.  The definition is correct because $\arctanh_p(x)\to\infty$ as $x\to1$ for any $p>0$.

Different variations of the functions $\sin_p$ and $\cos_p$ can be traced back to the 1879 paper of Lundberg \cite{Lind-Peetre2004,Lundberg} which remained forgotten until Jaak Peetre found it in 1995 \cite{Lindqvist95}. The next time  related functions seem to appear in a paper by David Shelupsky \cite{Shelupsky} and some of their values were computed by Burgoyne in 1964 \cite{Burgoyne}.  The generalized hyperbolic functions have been encountered by Peetre in \cite{Peetre72} in connection with the study of $K$-functional.  Later on, the function  $\sin_p$ was found to be the eigenfunction of a boundary value problem for one-dimensional $p$-Laplacian \cite{PEM89,DrMan99,Elbert79,Lindqvist95,Otani}.

More recently, monotonicity and convexity properties of the generalized  trigonometric functions and various inequalities for these functions have been  extensively studied by many authors. See \cite{BBK-Aequat,BBV,BV-JAT,BushEdm-Rocky,EdmGurkaLang-JAT,Qi-JAT13,KVZ,Wang-Chu,YinHuang} and references in these papers.
In particular, the monotonicity of $\sin_p(\pi_px)$ as a function of $p$ was required in \cite{BushEdm-Rocky} to fill a gap in the proof of basis properties of these functions given in \cite{BBCDG}.   Baricz, Bhayo and Vuorinen \cite{BBV} investigated convexity properties of  the functions $p\to\arcsin_p$, $p\to\arctan_p$ and their hyperbolic analogues and two-parameter generalizations.   They also proposed  the following conjecture.

\textbf{Conjecture (\cite{BBV})}.  The following Tur\'{a}n type
inequalities hold for all $p>2$ and $y\in(0,1)$:
$$
[\sin_p(y)]^2>\sin_{p-1}(y)\sin_{p+1}(y),
$$$$
[\cos_p(y)]^2>\cos_{p-1}(y)\cos_{p+1}(y),
$$$$
[\tan_{p}(y)]^2<\tan_{p-1}(y)\tan_{p+1}(y),
$$$$
[\sinh_p(y)]^2<\sinh_{p-1}(y)\sinh_{p+1}(y),
$$$$
[\tanh_p(y)]^2>\tanh_{p-1}(y)\tanh_{p+1}(y).
$$
The domain $(0,1)$ in this conjecture requires some explanation.  As we have seen above the functions $\sin_p(y)$, $\cos_p(y)$ and $\tan_{p}(y)$ are defined on $[0,\pi_p/2]$ as the inverse functions of the corresponding integrals. Note that for $y\in(\pi_p/2,\pi_p]$ the functions $\cos_p(y)$ and $\tan_{p}(y)$ are negative and their logarithmic convexity properties are not well defined. For the purposes of present investigation we will therefore restrict our attention to the intersection of the domains $[0,\pi_p/2]$ which is precisely $[0,1]$ since $\pi_{\infty}/2=1$.  On the other hand, the hyperbolic functions are defined as the inverse functions of the appropriate integrals for all $y\in\R_+$ so it is natural to extend the part of the above conjecture pertaining to hyperbolic functions to $y\in\R_+$.

In this paper we will prove that $p\to\sin_p(y)$ is logarithmically concave on $(0,\infty)$ for each $y\in[0,1]$, $p\to\sinh_p(y)$ is logarithmically convex and $p\to\tanh_p(y)$ is concave on $(0,\infty)$ for each  $y\in\R_+$.
These results confirm, strengthen and extend the above conjecture for $\sin_p$, $\sinh_p$ and $\tanh_p$.  Since $\cosh_p$ is a ratio of $\sinh_p$ and $\tanh_p$ we also conclude that $p\to\cosh_p$ is logarithmically convex. The conjecture for $\tan_p$ and $\cos_p$ will be demonstrated under stronger restriction $y\in(0,\log2)$.  We believe, however,  that it actually holds for all $y\in(0,1)$ but were unable to come up with a proof.

\paragraph{2. Auxiliary results.}
The following lemma will be our key tool for the forgoing investigation of the convexity properties of generalized trigonometric functions.
\begin{lemma}\label{lm:derivatives}
Suppose $I,J$ are finite or infinite open or closed subintervals
of $\mathbb{R}$.  Suppose $f(p,x)\in{C^2(J\times{I})}$ is strictly
monotone on $I$ for each fixed $p\in{J}$ so that $y\to
g(p,y):=f^{-1}(p,y)$ is well defined and monotone on $f(I)$ for
each fixed $p\in{J}$.  Then the following relations hold true
\begin{align}
&\frac{\partial}{\partial{p}}g(p,y)=-f'_{p}/f'_x,\label{eq:dgdalpha}
\\[5pt]
&\frac{\partial^2}{\partial{p^2}}g(p,y)=\left(f'_x\right)^{-2}
\left\{2f'_{p}f''_{xp}-f'_xf''_{pp}-\left(f'_{p}\right)^2f''_{xx}/f'_x\right\},
\label{eq:d2gd2alpha}
\\[5pt]
&\frac{\partial^2}{\partial{p^2}}\log{g(p,y)}=(xf'_x)^{-2}\left\{2xf'_{p}f''_{p{x}}-xf'_xf''_{pp}-x(f'_{p})^2f''_{xx}/f'_x-(f'_{p})^2\right\},
\label{eq:d2loggd2alpha}
\end{align}
where $x$ on the right is related to $y$ on the left by $y=f(p,x)$ or $x=g(p,y)$.
\end{lemma}
\textbf{Remark.} Formulas (\ref{eq:d2gd2alpha}) and (\ref{eq:d2loggd2alpha}) can also be written in the following form:
$$
\frac{\partial^2}{\partial{p^2}}g(p,y)=\frac{1}{2}\frac{\partial}{\partial{x}}\left(\frac{f'_{p}}{f'_{x}}\right)^{\!\!2}-\frac{\partial}{\partial{p}}\left(\frac{f'_{p}}{f'_{x}}\right)
$$
and
$$
\frac{\partial^2}{\partial{p^2}}\log{g(p,y)}=\frac{1}{2x}\frac{\partial}{\partial{x}}\left(\frac{f'_{p}}{f'_{x}}\right)^{\!\!2}
-\frac{1}{x}\frac{\partial}{\partial{p}}\left(\frac{f'_{p}}{f'_{x}}\right)
-\left(\frac{f'_{p}}{xf'_x}\right)^2.
$$
\textbf{Proof.} By definition of the inverse function we have:
\begin{equation}\label{eq:inverse}
f(p,g(p,y))=y
\end{equation}
Differentiating (\ref{eq:inverse}) with respect to $p$ while holding $y$ fixed we get:
$$
\frac{df}{dp}=\frac{\partial{f}}{\partial{p}}+\frac{\partial{f}}{\partial{x}}\frac{\partial{g}}{\partial{p}}=0~\Leftrightarrow~
\frac{\partial{g}}{\partial{p}}=-\frac{\partial{f}}{\partial{p}}\left(\frac{\partial{f}}{\partial{x}}\right)^{-1}.
$$
This proves (\ref{eq:dgdalpha}).  Further, differentiating $\partial{g}/\partial{p}$ once more with respect to $p$ yields:
\begin{multline*}
\frac{\partial^2{g}}{\partial{p}^2}=
-\left(\frac{\partial{f}}{\partial{x}}\right)^{-1}\frac{\partial}{\partial{p}}\left(\frac{\partial{f}}{\partial{p}}\right)
-\frac{\partial{f}}{\partial{p}}\frac{\partial}{\partial{p}}\left(\frac{\partial{f}}{\partial{x}}\right)^{-1}
=
\\
=-\left(\frac{\partial{f}}{\partial{x}}\right)^{-1}\left(
\frac{\partial^2{f}}{\partial{p}^2}+\frac{\partial^2{f}}{\partial{x}\partial{p}}\frac{\partial{g}}{\partial{p}}\right)
+\frac{\partial{f}}{\partial{p}}\left(\frac{\partial{f}}{\partial{x}}\right)^{-2}
\frac{\partial}{\partial{p}}\frac{\partial{f}}{\partial{x}}
\\
=-\left(\frac{\partial{f}}{\partial{x}}\right)^{-1}\left(
\frac{\partial^2{f}}{\partial{p}^2}+\frac{\partial^2{f}}{\partial{x}\partial{p}}\frac{\partial{g}}{\partial{p}}\right)
+\frac{\partial{f}}{\partial{p}}\left(\frac{\partial{f}}{\partial{x}}\right)^{-2}
\left(\frac{\partial^2{f}}{\partial{p}\partial{x}}+\frac{\partial^2{f}}{\partial{x}^2}\frac{\partial{g}}{\partial{p}}\right).
\end{multline*}
Substituting (\ref{eq:dgdalpha}) for $\partial{g}/\partial{p}$  into the last formula we obtain:
\begin{multline*}
\frac{\partial^2{g}}{\partial{p}^2}=
-\left(\frac{\partial{f}}{\partial{x}}\right)^{-1}\left(
\frac{\partial^2{f}}{\partial{p}^2}+\frac{\partial^2{f}}{\partial{x}\partial{p}}\left[-\frac{\partial{f}}{\partial{p}}\left(\frac{\partial{f}}{\partial{x}}\right)^{-1}\right]\right)
\\
+\frac{\partial{f}}{\partial{p}}\left(\frac{\partial{f}}{\partial{x}}\right)^{-2}
\left(\frac{\partial^2{f}}{\partial{p}\partial{x}}+\frac{\partial^2{f}}{\partial{x}^2}\left[-\frac{\partial{f}}{\partial{p}}\left(\frac{\partial{f}}{\partial{x}}\right)^{-1}\right]\right)
\\
=-\left(\frac{\partial{f}}{\partial{x}}\right)^{-1}\frac{\partial^2{f}}{\partial{p}^2}+2\left(\frac{\partial{f}}{\partial{x}}\right)^{-2}\frac{\partial{f}}{\partial{p}}\frac{\partial^2{f}}{\partial{x}\partial{p}}
-\left(\frac{\partial{f}}{\partial{x}}\right)^{-3}\frac{\partial^2{f}}{\partial{x}^2}\left(\frac{\partial{f}}{\partial{p}}\right)^2
\\
=-\left(\frac{\partial{f}}{\partial{x}}\right)^{-2}
\left[
\frac{\partial{f}}{\partial{x}}\frac{\partial^2{f}}{\partial{p}^2}
-2\frac{\partial{f}}{\partial{p}}\frac{\partial^2{f}}{\partial{x}\partial{p}}
+\left(\frac{\partial{f}}{\partial{x}}\right)^{-1}\frac{\partial^2{f}}{\partial{x}^2}\left(\frac{\partial{f}}{\partial{p}}\right)^2
\right]
\\
=\frac{1}{2}\frac{\partial}{\partial{x}}\left(\frac{f'_{p}}{f'_{x}}\right)^{\!\!2}-\frac{\partial}{\partial{p}}\left(\frac{f'_{p}}{f'_{x}}\right).
\end{multline*}
To prove (\ref{eq:d2loggd2alpha}) compute
$$
\frac{\partial^2}{\partial{p^2}}\log{g(p,y)}=\frac{g''_{pp}g-(g'_{p})^2}{g^2}.
$$
Substituting $g(p,y)=x$ and formulas for $g'_{p}$ and $g''_{pp}$ just derived into the above formula we obtain:
$$
\frac{\partial^2}{\partial{p^2}}\log{g(p,y)}
=\frac{1}{2x}\frac{\partial}{\partial{x}}\left(\frac{f'_{p}}{f'_{x}}\right)^{\!\!2}
-\frac{1}{x}\frac{\partial}{\partial{p}}\left(\frac{f'_{p}}{f'_{x}}\right)
-\left(\frac{f'_{p}}{xf'_x}\right)^2.
$$
Performing the differentiation we get (\ref{eq:d2loggd2alpha}). \hfill $\square$

We proceed with some standard definitions.  A positive function $f$ defined on a finite or infinite interval $I$ is said to be logarithmically convex, or log-convex, if its logarithm is convex, or equivalently,
$$
f(\lambda{x}+(1-\lambda)y)\leq[f(x)]^{\lambda}[f(y)]^{1-\lambda}~~\text{for all}~x,y\in{I}~~\text{and}~\lambda\in[0,1].
$$
The function $f$ is log-concave if the above inequality is reversed.  If the inequality sign is strict for $\lambda\in(0,1)$ then the appropriate property is called strict.  It is relatively straightforward to see from these definitions that log-convexity implies convexity and while log-concavity is implied by concavity but not vice versa.

The following corollaries specialize the formulas from
Lemma~\ref{lm:derivatives} to those generalize trigonometric
functions we will deal with in the present paper. Suppose $b>a>0$, where $b$ may equal $\infty$, are fixed.
If $p$ varies over $(a,b)$ the common domain of definition for the
families $\{\sin_p(y)\}_{p\in(a,b)}$, $\{\cos_p(y)\}_{p\in(a,b)}$ and $\{\tan_p(y)\}_{p\in(a,b)}$
is $[0,\pi_b/2]$, where $\pi_b=\infty$ for $0<b\le1$ and $\pi_{\infty}=2$ (we restrict our attention to the primary definitions as the inverse functions of the corresponding integrals).  The families $\{\sinh_p(y)\}_{p\in(a,b)}$, $\{\cosh_p(y)\}_{p\in(a,b)}$ and $\{\tanh_p(y)\}_{p\in(a,b)}$ are all defined on $[0,\infty)$.

\begin{corollary} \label{cr:sinp}
The function $p\to\sin_p(y)$ is log-concave on the interval $(a,b)$ for some
$y\in[0,\pi_b/2]$ iff for all $p\in(a,b)$
\begin{equation}\label{eq:lcsinp}
\frac{1}{x}\phi(p,x)^{p-1}\left(\int_0^x\phi'_p(p,t)dt\right)^{\!\!2}-2[\log\phi(p,x)]'_p\int_0^x\phi'_p(p,t)dt+\int_0^x\phi''_{pp}(p,t)dt\ge0,
\end{equation}
where $x=\sin_p(y)$.  It is concave on $(a,b)$ for some
$y\in[0,\pi_b/2]$ iff for all $p\in(a,b)$
\begin{equation}\label{eq:csinp}
\frac{x^p}{x}\phi(p,x)^{p-1}\left(\int_0^x\phi'_p(p,t)dt\right)^{\!\!2}-2[\log\phi(p,x)]'_p\int_0^x\phi'_p(p,t)dt+\int_0^x\phi''_{pp}(p,t)dt\ge0.
\end{equation}
Here  $\phi(p,t)=(1-t^p)^{-1/p}$ and
$$
\frac{\phi'_p(p,t)}{\phi(p,t)}\!=\!\frac{1}{p^2}\log(1-t^p)+\frac{t^p\log{t}}{p(1-t^p)},~~
\phi''_{pp}(p,t)\!=\!\frac{\phi'_p(p,t)^2}{\phi(p,t)}-\frac{2}{p}\phi'_p(p,t)+\phi(p,t)\frac{t^p\log^2{t}}{p(1-t^p)^2}.
$$
The corresponding property is strict if and only if the inequality sign is strict.
\end{corollary}

\textbf{Proof.}  Write $g(p,y)=\sin_p(y)$. A necessary and
sufficient condition for log-concavity of the smooth function $p\to{g(p,y)}$
is $[\log(g)]''_{pp}\leq0$.  To compute $[\log(g)]''_{pp}$
substitute $f(p,x)=F_p(x)$ defined in (\ref{eq:arcsin}) into
(\ref{eq:d2loggd2alpha}) and notice that
$$
f'_x(p,x)=(1-x^p)^{-1/p}=\phi(p,x),~~
f''_{xx}(p,x)=(1-x^p)^{-1-1/p}x^{p-1}=\frac{\phi(p,x)x^{p-1}}{1-x^p},
$$
$$
\frac{x}{f'_x}f''_{xx}+1=\frac{x(1-x^p)^{-1-1/p}x^{p-1}}{(1-x^p)^{-1/p}}+1=\frac{x^p}{1-x^p}+1=\frac{1}{1-x^p}=\phi(p,x)^p,
$$
$$
f''_{xp}=f''_{px}=\phi'_p(p,x),~~
f'_p=\int\limits_{0}^{x}\phi'_p(t)dt,~~
f''_{pp}=\int\limits_{0}^{x}\phi''_{pp}(t)dt.
$$
Hence $[\log(g)]''_{pp}\leq0$ reduces to (\ref{eq:lcsinp}).
Similarly using (\ref{eq:d2gd2alpha}) we get (\ref{eq:csinp}).
Formulas for derivatives are obtained by straightforward
differentiation.  \hfill$\square$

In precisely the same manner we can derive the next three
corollaries whose proofs are omitted.

\begin{corollary} \label{cr:ptan}
 The function
$p\to\tan_p(y)$ is log-convex on the interval $(a,b)$ for some
$y\in[0,\pi_b/2]$ iff for all $p\in(a,b)$
\begin{equation}\label{eq:lctanp}
(1/x-(p-1)x^{p-1})\left(\int_0^x\theta'_p(p,t)dt\right)^{\!\!2}+\frac{2x^p\log{x}}{1+x^p}\int_0^x\theta'_p(p,t)dt+\int_0^x\theta''_{pp}(p,t)dt\le0,
\end{equation}
where $x=\tan_p(y)$ and $\theta(p,x)=(1+x^p)^{-1}$.  It is convex
on $(a,b)$ for some $y\in[0,\pi_b/2]$ iff for all $p\in(a,b)$
\begin{equation}\label{eq:ctanp}
px^{p-1}\left(\int_0^x\theta'_p(p,t)dt\right)^{\!\!2}+\frac{2x^p\log{x}}{1+x^p}\int_0^x\theta'_p(p,t)dt+\int_0^x\theta''_{pp}(p,t)dt\le0.
\end{equation}
Here
$$
\theta'_p(p,t)=\frac{-t^p\log{t}}{(1+t^p)^2},~~~~\theta''_{pp}(p,t)=\frac{-t^p(1-t^p)\log^2{t}}{(1+t^p)^3}.
$$
The corresponding property is strict if and only if the inequality sign is strict.
\end{corollary}
\begin{corollary} \label{cr:psinhp}
The function $p\to\sinh_p(y)$ is log-convex on the interval $(a,b)$ for
some $y\in[0,+\infty)$ iff for all $p\in(a,b)$
\begin{equation}\label{sinh}
\frac{1}{1+x^p}\left(\int_0^x\lambda'_p(p,t)dt\right)^2-2x\lambda'_p(p,x)\int_0^x\lambda'_p(p,t)dt
+\frac{x}{(1+x^p)^{1/p}}\int_0^x\lambda''_{pp}(p,t)dt\le0.
\end{equation}
Here $x=\sinh_p(y),$ $\lambda(p,t)=(1+t^p)^{-1/p}$ and
$$
\frac{\lambda'_p(p,t)}{\lambda(p,t)}\!=\!\frac{\log(1+t^p)}{p^2}-\frac{t^p\log{t}}{p(1+t^p)},
~~~
\lambda''_{pp}(p,t)\!=\!\frac{\lambda'_p(p,t)^2}{\lambda(p,t)}-\frac{2}{p}\lambda'_p(p,t)
-\lambda(p,t)\frac{t^{p}\log^2{t}}{p(1+t^p)^2}.
$$
The corresponding property is strict if and only if the inequality sign is strict.
\end{corollary}

\begin{corollary} \label{cr:tanhq}
The function $p\to\tanh_p(y)$ is concave on the interval $(a,b)$
for some $y\in[0,\infty)$ iff for all $p\in(a,b)$
\begin{equation}\label{eq:qsinpq}
\frac{px^{p-1}}{(1-x^p)}\left(\int\limits_0^x\alpha'_p(p,t)dt\right)^2
-\frac{2x^p\log{x}}{(1-x^p)^{2}}\int\limits_0^x\alpha'_p(p,t)dt+\frac{1}{(1-x^p)}\int\limits_0^x\alpha''_{pp}(p,t)dt\geq0,
\end{equation}
where $x=\tanh_p(y).$ Here
$$\alpha(p,t)=\frac{1}{1-t^p},~~~\alpha'_p(p,t)=\frac{t^p\log{t}}{(1-t^p)^{2}},
~~~\alpha''_{pp}(p,t)=\frac{t^p(t^p+1)\log^2{t}}{(1-t^p)^{3}}.
$$
The corresponding property is strict if and only if the inequality sign is strict.
\end{corollary}

The next lemma is a guise of the  monotone L'H\^{o}spital rule \cite{AVV,Pinelis}. As before we allow the value $b=\infty$.
\begin{lemma}\label{lm:LHospital}
Suppose $u$, $v$ are continuously differentiable functions defined on a real interval $(a,b)$, $u(a)=v(a)=0$ and $vv'>0$ on $(a,b)$. If $u'/v'$ is decreasing on $(a,b)$ then $u/v>u'/v'$ on $(a,b)$.
\end{lemma}
\textbf{Proof.} According to monotone L'H\^{o}spital rule decrease
of $u'/v'$ implies that $u/v$ is also decreasing.  On the other
hand, by the quotient rule (see also \cite[formula
(1.1)]{Pinelis})
$$
v^2\left(\frac{u}{v}\right)'=\left(\frac{u'}{v'}-\frac{u}{v}\right)vv'
$$
so that the expression in parentheses must be negative. \hfill
$\square$

We will also need the following estimate.
\begin{lemma}\label{lm:tan}
Suppose $p>1$, $0<s<1$.  Then
\begin{equation}\label{in:tan}
\frac{sp^3}{(p+1)^2}\left(\frac{1}{(p+1)^2}-\frac{\log^2{s}}{p^2}\right)<\int_0^1u^{1/p}\frac{(1-su)}{(1+su)^3}\log^2(su)du.
\end{equation}
\end{lemma}
\textbf{Proof.} Define
$$
t(s)=s\left(\frac{1}{(p+1)^2}-\frac{\log^2{s}}{p^2}\right).
$$
Taking derivative yields
$$
t'_s=\frac{1}{(p+1)^2}-\frac{\log^2{s}}{p^2}-\frac{2\log{s}}{p^2}.
$$
The equation  $t'_s=0$ has two roots  $s_{1,2}=\exp\{-1\pm\sqrt{1+p^2/(p+1)^2}\}$.
Hence, one of the roots lies in $(0,1)$ and $t(s)$ is decreasing on $(0,s_*)$ and increasing on $(s_*,1)$ for some $0<s_*<1$.  Since $t(0)=0$ and $t(1)>0$ the maximum is attained at the point $s=1$.  On the other hand, the function  $s\to(1-su)\ln^2(su)/(1+su)^3$  decreases on $[0,1]$ and attains its minimum at the point $s=1$. This implies that we only need to prove the inequality
\begin{equation}\label{in:tan1}
\frac{p^3}{(p+1)^4}<\int_0^1u^{1/p}\frac{(1-u)}{(1+u)^3}\log^2{u}du.
\end{equation}
In view of
 $$
 \int_0^1u^{1/p}\log^2{u}du=\frac{2p^3}{(1+p)^3}
 $$
 we can rewrite (\ref{in:tan1}) as
\begin{equation}\label{in:tan2}
\int_0^1
u^{1/p}\log^2{u}\left(\frac{(1-u)}{(1+u)^3}-\frac{1}{2(p+1)}\right)du>0.
\end{equation}
The derivative of the integrand with respect to $p$ equals
$$
\frac{u^{1/p}\log^2{u}}{p^2}\left\{\frac{p^2}{2(p+1)^2}-\log{u}\left(\frac{(1-u)}{(1+u)^3}-\frac{1}{2(p+1)}\right)\right\}.
$$
The expression in braces increases in $p$ for each fixed $u\in(0,1]$ and is easily seen to be positive for $p=1$.
Therefore, it is positive for all $p\ge1$ and $u\in(0,1]$. This shows that the left hand side of (\ref{in:tan2}) increases in $p$ for $p\ge1$. When $p=1$ computation gives
$$
\int_0^1u^{1/p}\log^2{u}\left(\frac{(1-u)}{(1+u)^3}-\frac{1}{2(p+1)}\right)du=\frac{\pi^2}{3}-\log
4-\frac{3}{2}\zeta(3)-\frac{1}{16}>0,
$$
where $\zeta$ is Riemann's zeta function.  The last inequality completes the proof.\hfill $\square$

\bigskip

\paragraph{Main results.}

Our first theorem is concerned with the generalized sine.
\begin{theorem}\label{th:sinp}
For each fixed $y\in(0,1)$ the function $p\to\sin_p(y)$ is strictly log-concave on $(0,\infty)$.
\end{theorem}
\textbf{Proof.} We need to prove (\ref{eq:lcsinp}) with strict
sign. Consider the following quadratic
$$
F(z)=\frac{\phi(p,x)^{p-1}}{x}z^2-2\eta(p,x)z+\int_0^x\phi''_{pp}(p,t)dt,
$$
where $\eta(p,x)=[\log\phi(p,x)]'_p$.  Explicit expressions for $\phi$, $\eta$ and $\phi''_{pp}$ is given in corollary
\ref{cr:sinp}.  We need to show that $F(z)>0$ for $z=\int_0^x\phi'_{p}(p,t)dt$. We will show that in fact this
inequality holds for all real $z$. Indeed, $\phi(p,x)^{p-1}/x>0$ and it remains to prove that
$$
\frac{D}{4}=\eta^2-\frac{\phi(p,x)^{p-1}}{x}\int_0^x\phi''_{pp}(p,t)dt<0~\Leftrightarrow~G(x):=\frac{x\eta^2}{\phi(p,x)^{p-1}}-\int_0^x\phi''_{pp}(p,t)dt<0.
$$
Here $D$ denotes the discriminant of $F(z)$.  Clearly, $G(0)=0$.
Further, elementary but long computation reveals
\begin{multline*}
G'(x)=\frac{2\eta(p,x)(1-x^p+px^p\log{x})}{p(1-x^p)^{1+1/p}}-\frac{px^p\eta(p,x)^2}{(1-x^p)^{1/p}}-\frac{x^p\log^2{x}}{p(1-x^p)^{2+1/p}}
\\
=(1-x^p)^{-1/p}\left\{\frac{2}{p}\eta(p,x)-px^p\eta(p,x)^2+\frac{2\eta(p,x)x^p\log{x}}{1-x^p}-\frac{x^p\log^2{x}}{p(1-x^p)^{2}}\right\}.
\end{multline*}
Substituting
$$
\eta(p,x)=\frac{1}{p^2}\log(1-x^p)+\frac{x^p\log{x}}{p(1-x^p)}
$$
into the above formula  we get after some more algebra:
$$
G'(x)=(1-x^p)^{-1/p}\left\{\frac{2x^p\log{x}}{p^2(1-x^p)}+\frac{2}{p^3}\log(1-x^p)-\frac{x^p}{p}(\log(1-x^p)^{1/p}-\log{x})^2\right\}.
$$
Each term in braces is negative so that $G'(x)<0$ which implies that $G(x)<0$ completing the proof of the theorem.\hfill $\square$

Extensive numerical evidence supports the following assertion.

\textbf{Conjecture.}  There exists $p_0\in(0,1)$ such that the
function $p\to\sin_p(y)$ is strictly concave on $(p_0,\infty)$ for
all $y\in(0,1)$.  If $p\in(0,p_0)$ concavity is violated for some
$y\in(0,1)$.

\bigskip
\begin{theorem}\label{th:tan}
For each fixed $y\in(0,\log{2})$ the function $p\to \tan_p(y)$ is
strictly log-convex on $(1,\infty)$.
\end{theorem}

\textbf{Proof.} First we show that the claim of the Theorem is equivalent to validity of inequality
(\ref{eq:lctanp}) with strict sign for $x\in(0,1)$ and $p>1$.  Indeed, $y\to\tan_p(y)$ is increasing so that $\tan_1(y)<\tan_1(\log{2})=1$ for $y\in(0,\log{2})$.  Further, $\tan_p(y)<\tan_1(y)<1$ because
$$
\int_0^{\tan_p(y)}\frac{dt}{1+t^p}=y=\int_0^{\tan_1(y)}\frac{dt}{1+t}<\int_0^{\tan_1(y)}\frac{dt}{1+t^p}.
$$

Since $1+(1-p)x^p<1$ it suffices to prove the inequality
\begin{equation}\label{in:tan5}
\left(\int_0^x\theta'_p(p,t){dt}\right)^2
+\frac{2x^{p+1}}{1+x^p}\log(x)\int_0^x\theta'_p(p,t){dt}+x\int_0^x\theta''_{pp}(p,t){dt}<0.
\end{equation}
Note that $\theta'_p(p,t)>0$, $\theta''_{pp}(p,t)<0$ for $0<t<1$.
Further, by replacing the denominator in $\theta'_p(p,t)$ by $1$ we get the estimate
\begin{equation}\label{eq:psi-est}
0<\int_0^x \theta'_p(p,t)dt<\psi:=-\int_0^xt^p\log(t)dt=\frac{x^{p+1}}{p+1}\left(\frac{1}{p+1}-\log(x)\right).
\end{equation}
Consider the following quadratic derived from (\ref{in:tan5}):
$$
F(w)=w^2+bw+c,~~\text{where}~b=\frac{2x^{p+1}\log(x)}{1+x^p},~~c=x\int_0^x\theta''_{pp}(p,t)dt.
$$
We need to show that $F(w)<0$ for $w=\int_0^x\theta'_p(p,t)dt$.
Since $c$ is negative $F(0)<0$ and in view of (\ref{eq:psi-est})
it suffices to demonstrate that $F(\psi)<0$ or
\begin{equation}\label{in:tan6}
-c>\psi^2+\psi b.
\end{equation}
Denote  $s=x^p$.  Then after some rearrangement we get
\begin{equation}\label{in:tan7}
\psi^2+\psi{b}=\frac{s^2x^2}{(p+1)^2}\left(\frac{-\log(s)}{p}+\frac{1}{p+1}\right)
\left(\frac{-\log(s)}{p}+\frac{1}{p+1}+\frac{2(p+1)\log(s)}{p(s+1)}\right).
\end{equation}
Performing change of variable $t=ux$ in the integral representing $c$ we obtain
$$
-c=x\int_0^x \frac{t^p (1-t^p)\log^2(t)dt}{(1+t^p)^3}
=x\int_0^1\frac{\log^2(u^px^p)pu^{p-1}x^p(1-x^pu^p)xudu}{p^3(1+x^pu^p)^3}.
$$
By writing $s=x^p$, $z=u^p$ the last integral transforms into
\begin{equation}\label{in:tan8}
-c=\frac{x^2s}{p^3}\int_0^1z^{1/p}\frac{(1-sz)}{(1+sz)^3}\log^2(sz)dz.
\end{equation}
By dropping some terms in (\ref{in:tan7}) we have
$$
\psi^2+\psi{b}<\frac{s^2x^2}{(p+1)^2}\left(\frac{1}{(p+1)^2}-\frac{\log^2(s)}{p^2}\right).
$$
The right hand side of this inequality does not exceed $-c$ by Lemma~\ref{lm:tan}.\hfill $\square$

We remark here that log-convexity of $p\to\tan_p(y)$ does not hold for all $p>0$ even if $y$ is restricted to $(0,\log2)$ as evidenced by numerical experiments.

\begin{theorem} \label{th:cos}
For each fixed $y\in(0,\log 2)$ the function $p\to \cos_p(y)$ is
strictly log-concave on $(1,\infty)$.
\end{theorem}

\textbf{Proof.} By Theorem~\ref{th:tan} the function $\tan_p(y)$ is log-convex on $(1,\infty)$, i.e.
$(\log\tan_p(y))''_{pp}>0$. On the other hand, according to Theorem~\ref{th:sinp} $(\log\sin_p(y))''_{pp}<0$ and the claim now follows from the formula (see (\ref{eq:tanp}))
$$
(\log\cos_p(y))''_{pp}=(\log\sin_p(y))''_{pp}-(\log\tan_p(y))''_{pp}.~~~~~~~~~~~\square
$$

\begin{theorem}\label{th:shq}
For each fixed $y\in(0,\infty)$ the function $p\to \sinh_p(y)$ is strictly  log-convex on $(0,\infty)$.
\end{theorem}
\textbf{Proof.}  By Corollary~\ref{cr:psinhp} we need to prove (\ref{sinh}).
After dividing (\ref{sinh}) by $-x\lambda(p,x)<0$ we get
\begin{equation}\label{in:sinh}
v(p,x)\mu(p,x)^2+2w(p,x)\mu(p,x)-\int_0^x\lambda''_{pp}(p,t)dt>0,
\end{equation}
where
$$
v(p,x)=-\frac{(1+x^p)^{-1+1/p}}{x},
~~~w(p,x)=\frac{1}{p^2}\log(1+x^p)-\frac{x^p\log(x)}{p(1+x^p)},
~~~\mu(p,x)=\int_0^x \lambda'_p(p,t)dt.
$$
We denote the left hand side of (\ref{in:sinh}) by $H(p,x)$ and abbreviate $v=v(p,x)$, $\mu=\mu(p,x)$,  $w=w(p,x)$. Differentiation gives
$$
H'_x(p,x)=v'_x\mu^2+2(v\lambda'_p+w'_x)\mu+(2w\lambda'_p-\lambda''_{pp}).
$$
We will demonstrate that in fact $H'_x(x,p)>0$ for all real $\mu$. Indeed,
$$
v'_x=\frac{(1+x^p)^{-2+1/p}(1+px^p)}{x^2}>0
$$
and it remains to prove that
$$
\frac{D}{4}=(v\lambda'_p+w'_x)^2-(2w\lambda'_p-\lambda''_{pp})v'_x<0,
$$
where $D$ denotes the discriminant of the quadratic $H'_x(x,p)$ viewed as a function of $\mu$.  Straightforward calculations reveal:
$$
v=-\frac{1}{\lambda(1+x^p)x},~~\lambda_p'=\lambda w,~~w'_x=\frac{-x^{p-1}\log x}{(1+x^p)^2},
$$
$$
\lambda_{pp}''=\lambda\left(w^2-\frac{2}{p}w-\frac{x^p\log^2(x)}{p(1+x^p)^2}\right),
~~v'_x=\frac{1+px^p}{\lambda{x^2}(1+x^p)^2}.
$$
Substitution yields
\begin{equation}{\label{eq:discriminant}}
\frac{D}{4}=\frac{1}{x^2(1+x^p)^2}\left\{-px^pw^2-\frac{2w}{p}\left(1+px^p-\frac{x^p\log(x^p)}{1+x^p}\right)
-\frac{x^p\log^2(x)}{p(1+x^p)^2}\right\}.
\end{equation}
We only need to show that the middle term in braces is positive or
$$
\frac{2w}{p}\left(1+px^p-\frac{x^p\log(x^p)}{1+x^p}\right)=\frac{2\left((1+z)\log(1+z)-z\log z\right)
\left(1+pz-\frac{z\log z}{1+z}\right)}{p^3(1+z)}>0,
$$
where $z=x^p$. Indeed, this follows from  the inequalities
$$
(1+z)\log (1+z)-z\log z>0,
$$
$$
1+pz-\frac{z\log z}{1+z}>1+z-\frac{z\log(z)}{1+z}>0
$$
valid for any $z>0$, $p>1$.  Hence, $H'_x(p,x)>0$ implying $H(p,x)>H(p,0)=0$ which completes the proof of the theorem.\hfill$\square$

\begin{theorem} \label{th:tanh}
For each fixed $y\in(0,\infty)$ the function $p\to \tanh_p(y)$ is
strictly  concave on $(0,\infty)$.
\end{theorem}

\textbf{Proof.} In order to prove concavity of $p\to \tanh_p(y)$
we need to show the validity of (\ref{eq:qsinpq}). The trick with quadratic
used in previous theorems fails here, so resort to another proof. The first term in (\ref{eq:qsinpq}) is positive so if we drop the expression on the left becomes smaller. Hence, if we can prove it is still positive we are done.
This amounts to showing that
$$
\frac{u}{v}:=\left(\int_0^x\frac{t^p(t^p+1)[\log(1/t)]^2}{(1-t^p)^{3}}dt\right)\!\Big/\!
\left(\int_0^x\frac{t^p\log(1/t)dt}{(1-t^p)^{2}}\right)>\frac{2x^p\log(1/x)}{(1-x^p)}.
$$
We have
$$
\frac{u'}{v'}=\frac{(x^p+1)\log(1/x)}{1-x^p}.
$$
It is easy to check by taking derivative that this function
decreases on $(0,1)$.  Clearly $u(0)=v(0)=0$ and $vv'>0$, so that
we are in the position to apply Lemma~\ref{lm:LHospital} yielding
$$
\left(\int_0^x\frac{t^p(t^p+1)[\log(t)]^2}{(1-t^p)^{3}}dt\right)\Big/\left(\int_0^x\frac{t^p\log(t)dt}{(1-t^p)^{2}}\right)>
\frac{(x^p+1)\log(1/x)}{1-x^p}>\frac{2x^p\log(1/x)}{(1-x^p)}.~~\square
$$

\begin{theorem} \label{th:cosh}
For each fixed $y\in(0,\infty)$ the function $p\to \cosh_p(y)$ is
strictly  log-convex on $(0,\infty)$.
\end{theorem}

\textbf{Proof.} By Theorem~\ref{th:tanh} the function
$p\to\tanh_p(y)$ is concave on $(0,\infty)$ and thus also log-concave, i.e. $(\log\tanh_p(y))''_{pp}<0.$
On the other hand, according to Theorem~\ref{th:shq}  that $(\log\cosh_p(y))''_{pp}>0$.  The claim now follows from the formula $$
(\log\cosh_p(y))''_{pp}=(\log\sinh_p(y))''_{pp}-(\log\tanh_p(y))''_{pp}.~~~~~~~\square
$$

\paragraph{Acknowledgements.} We thank Sergei Sitnik for numerous useful discussion regarding the manuscript.  The research was supported by Far Eastern Branch of the Russian Academy of Sciences (project 12-II-CO-01M-002) and Russian Foundation for Basic Research (project 13-01-12404-ofi\_m2).

\end{document}